\documentstyle{article}
\def\matrix{\begin{array}{c}} \def\endmatrix{\end{array}}
\def\dfrac{\displaystyle\frac}
\def\dint{\displaystyle\int}
\title{ The harmonic product of $\delta (x_{1},\ldots, x_{n})$ and
$\delta (x_{1})$\\ and two combinatorial identities
\thanks{Project supported by the National Natural Science Foundation of
China}}
\author{ Ya-Qing Li}
\date{}
\begin{document}
\maketitle
\baselineskip 18pt
\begin{abstract}

In the framework of nonstandard analysis, Bang-He Li and the author
defined the product of any two distributions on $R^n$ via their harmonic
representations. The
product of $\delta (x_{1},\ldots, x_{n})$ and
$\delta (x_{1})$ was calculated by Kuribayashi and the author in [LK].
In this paper,
the result of [LK] is improved to
$$\delta (x_{1},\ldots, x_{n})\circ \delta (x_{1})
=\dfrac{1}{2\pi\rho} \delta (x_{1},\ldots, x_{n})\;\,\mbox{
mod}\;\;\mbox{infinitesimals}$$
where $\rho$ is a positive infinitesimal.
Moreover two combinatorial identities are obtained as byproducts. 


03H05 46F05 46F10 46F20
\end{abstract}

Thirty years ago, Bremermann and Durand [BD] defined the products of 
distributions with one variable by using analytic representations. It was 
shown by Itano [I1][I2] and further by Bang-He Li and the author [LL1] that
this multiplication is very broad, i.e. if the product of two distributions 
exists for several other multiplications, then the same product is obtained  
for this multiplication. So a problem that interested people was 
`` What is a generalization of this multiplication to distributions with 
several variables?".

Itano [I2] showed by an example that for distributions with more than 
one variables, analytic representation can not offer well-defined 
multiplication. Bang-He Li and the author [LL2] at last found that a suitable 
generalization of this multiplication is the one via harmonic representations. 
Because for 
one variable, harmonic and analytic representations are essentially the same.
Bang-He Li [L] adapted the multiplication of Bremermann and Durand into the 
framework of nonstandard analysis. Its generalization to multiple variables
via harmonic 
representations [LL2] was also written in this framework. The merit to use 
nonstandard analysis is that one needs not to worry about the existence of 
products anymore, and when taking the finite part (if exists), we return to 
some kind of standard product.

In the case of one variable, systematic results for singular distributions 
have been obtained by the author. We quote only the survey 
paper [LL3]. For multiple variables, only few calculations have been 
made in [LL3] and the paper of Kuribayashi and the author 
[LK]. In this paper, we improve the result of [LK] to the neatest form.

Denote by ${\cal D}(R^{n})$ the Schwartz spase consisting of complex-valued
$C^\infty$-functions on $R^n$ with compact supports, and ${\cal D'}(R^{n})$
its dual, i.e. the space of Schwarts distributions.

For $T\in {\cal D'}(R^{n})$, its harmonic representation $\hat{T}$ is a
harmonic function on $R^n\times R_+$ such that
$$\lim_{y\longrightarrow 0_+} \int_{R^n}\hat{T}(x_1,\cdots,x_n,y)
\phi(x_1,\cdots,x_n)\,dx_1\cdots dx_n
=\int_{R^n}T(x_1,\cdots, x_n)\phi(x_1,\cdots, x_n)\,dx_1\cdots dx_n\;$$
for any $\phi(x_1,\cdots,x_n)\in {\cal D}(R^{n})$.

Harmonic representation exists for any $T\in {\cal D'}(R^{n})$, and the
difference of two such representations extends to a harmonic function on
$R^n\times R$ skew-symmetric for $y$ (cf. [LL4]).

Denote by $C$ the complex field, $\:{}^{*}C$ a nonstandard model of $C, \,\,R$
the real field and $\rho \in {}^{*}R$ a positive infinitesimal.\\
Let
$$\array{cl} {}^{\rho}C & =\{ x \in\, {}^{*}C | \mbox{ for some finite 
integer}\;\; n,\; |x| < \rho^{-n}\}\\\\
 \theta & = \mbox{the set of all infinitesimals in}\, {}^{*}C \\\\
 {}^{\rho}C' & = {}^{\rho}C/\theta
 \endarray$$
then  ${}^{\rho}C'$ is a complex vector space, and we call a complex linear 
functional of ${\cal D}(R^{n}) \longrightarrow
{}^{\rho}C'$ a hyperdistribution on $R^{n}$.

Suppose $S,\, T \in {\cal D}'(R^{n}), \hat{S},\hat{T}$ are harmonic
representations of $S$ and $T$, and ${}^{*}{\hat{S}}$, ${}^{*}{\hat{T}}$ are
the nonstandard extensions of $\hat{S}$ and $\hat{T} $ respectively.

Let $\psi : {}^{\rho}C \longrightarrow {}^{\rho}C'$ be the homomorphism
modulo $\theta$.\, Then
$$ \varphi \longrightarrow \psi ( < \hat{S}(x,\rho)\hat{T}(x,\rho), \varphi(x) >)
$$
defines a complex linear functional of\, ${\cal D}(R^{n})\longrightarrow
{}^{\rho }C'$, i.e. a hyperdistribution on $R^n$, and we call this hyperdistribution
the harmonic product of $S$ and $T$, denoted by $S\circ T$ in [LL1]
(see also [O]).\\

Let $\delta (x_{1},\ldots, x_{n})$ be the $\delta-$function on $R^n$,
 The harmonic product of
 $\delta (x_{1},\ldots, x_{n})$ and $\delta(x_1)$
 has been calculated in [LK] to get
$$\delta (x_{1},\ldots, x_{n})\circ \delta (x_{1})
=A(1, n) \delta (x_{1},\ldots, x_{n})$$
where
$$A(1, n)=
\frac{2\pi}{c_1 c_{n}\rho}\prod_{j=1}^{n-3}\int_{0}^{\pi}\sin^j
\theta\,d\theta\int_0^{\infty}\int_0^{\pi}\frac{t^{n}\sin^{n-2}\xi\,dt\,d\xi}
{(1+t^2)^{\frac{n+1}{2}}(1+t^2\cos^{2}\xi)}$$
and
$$c_n=\frac{\pi^{\frac{n+1}{2}}}{\Gamma(\frac{n+1}{2})}$$
Furthermore, for odd $n$,
$$A(1,2k+1)=\frac{2k}{\rho\pi}\big(\frac{1}{4}+\sum_{j=1}^{k-1}\sum_{p=0}
^{j-1}\left(\matrix k-1\\ j\\\endmatrix\right)\left(\matrix j-1\\ p\\
\endmatrix\right)\frac{(-1)^{j}(2p+1)!!(2j-2p-1)!! }
{(2j+2)!!(2p+1)}\big)$$
and for even $n$,
$$\begin{array}{ll}
&A(1,2k+2)=\dfrac{2k+1}{2\rho\pi}
\big(1+\sum_{j=0}^{k-1}\sum_{r=0}
^{j}\sum_{p=0}^{j-r}\sum_{s=0}^{r}\sum_{h=0}^{p+s}\left(\matrix 2k\\
k+1+j\\ \endmatrix\right)
\left(\matrix j-r\\ p\\ \endmatrix\right)\left(\matrix r\\ s\\
\endmatrix\right)\left(\matrix p+s\\ h\\\endmatrix\right)
\\ \\&\qquad\quad\qquad\quad\times
\dfrac{(-1)^{j+p+r+1} \Gamma(k-j) \Gamma(1-h+(p+s+j+1)/2) }
{2^{2k-2+p+s-j}\Gamma(k-h+(p+s-j+3)/2)}\big)\end{array}$$
where $\Gamma(x)$ is the Gamma function.

It was also calculated in [LK] that
$$A(1,2k+1)=\dfrac{1}{2\pi\rho},\;\; \mbox {for} \; k=1,2,3,4,5,6$$
and
$$A(1,2k+2)=\dfrac{1}{2\pi\rho},\;\; \mbox {for} \; k=0,1,2,3$$
It is a result of Bang-He Li [L] that
$$A(1,1)=\frac{1}{2\pi\rho}, \;\;\mbox{i.e.}\;\, \delta(x_1)\circ\delta(x_1)
=\frac{1}{2\pi\rho}\delta(x_1)
$$
(see also [O] [DR])

So it might be conjectured that $A(1,n)=\frac{1}{2\pi\rho}$ for 
all $n\geq 1$. Here we find a very simple way to prove that it is indeed 
so, i.e. we have

{\bf Theorem 1}.\qquad For any $n\in N$,
$$\delta (x_{1},\ldots, x_{n})\circ \delta (x_{1})
=\dfrac{1}{2\pi\rho} \delta (x_{1},\ldots, x_{n})
$$

{Proof}. \qquad A harmonic representation of $\delta (x_{1},\ldots, x_{n})$
is the Poisson kernel
$$\hat {\delta} (x_{1},\ldots, x_{n}; y)
={c_n}^{-1}y(|x|^2+y^2)^{-\frac{n+1}{2}},\;\;
y>0$$
and a harmonic representation of $\delta (x_{1})$ is
$$\hat {\delta} (x_1, y)={c_1}^{-1}y (x_1^2+y^2)^{-1}$$
Let $\phi\in {\sl D}(R^n)$, it has been proved in [LK] that
$$\begin{array}{ll}
&\;\;\,\;\;\;\dint_{R^n}\hat {\delta} (x_{1},\ldots, x_{n}; \rho)\hat {\delta} (x_1, \rho)
\phi(x)\,dx\\
&=\phi(0)\dint_{R^n}\hat {\delta} (x_{1},\ldots, x_{n}; \rho)
\hat {\delta} (x_1, \rho)dx\;\;\;\mbox{mod}\;\;\mbox{infinitesimals}\end{array}
$$
So, for $n\geq 3$
$$\begin{array}{ll}A(1,n)
&=\dfrac{\rho^{2}}{c_{1}c_{n}}\int_{R^{n}}\dfrac{dx}
{(x_1^{2}+x_2^2+\cdots+x_{n}^2+\rho^{2})^{\frac{n+1}{2}}(x_{1}^{2}+\rho^2)}\\ \\
&=\dfrac{\rho^{2}}{c_{1}c_{n}}\int_{R^{n-2}}\dfrac{dx_{1}\cdots dx_{n-2}}
{(x_{1}^2+\rho^2)}\,\int_{R^2}\dfrac{dx_{n-1}dx_{n}}{(
x_1^{2}+x_2^2+\cdots+x_{n}^2+\rho^{2})^{\frac{n+1}{2}}}
\end{array}$$
By using polar coordinate, we have
$$\begin{array}{ll}&\dint_{R^2}\dfrac{dx_{n-1}dx_{n}}{(
x_1^{2}+x_2^2+\cdots+x_{n}^2+\rho^{2})^{\frac{n+1}{2}}}=
\int_{0}^{2\pi}\,d\theta\int_0^{\infty}\dfrac{r\,dr}{(x_{1}^2+x_{2}^2
+\cdots+x_{n-2}^2+\rho^{2}+r^2)^{\frac{n+1}{2}}}\\
=&\dfrac{2\pi}{(n-1)(x_{1}^2+x_{2}^2+\cdots+x_{n-2}^2+\rho^{2})^
{\frac{n-1}{2}}}\end{array}$$
Now
$$c_n=\frac{2\pi}{n-1} c_{n-2}$$
Hence

$$A(1,n)
=\dfrac{\rho^{2}}{c_{1}c_{n-2}}\int_{R^{n-2}}\dfrac{dx_1\cdots dx_{n-2}}
{(x_1^{2}+x_2^2+\cdots+x_{n-2}^2+\rho^{2})^{\frac{n-1}{2}}(x_{1}^{2}+\rho^2)}
=A(1,n-2) $$
and the proof is complete by using the known results
$$A(1,1)=A(1,2 )=\frac{1}{2\pi\rho}$$

As corollaries of Theorem 1 and the results of [LK], we obtain two 
combinatorial identities

{\bf Theorem 2}.\qquad For any $k\in N,\,k\geq 1$,
$$\sum_{j=1}^{k-1}\sum_{p=0}
^{j-1}\left(\matrix k-1\\ j\\\endmatrix\right)\left(\matrix j-1\\ p\\
\endmatrix\right)\frac{(-1)^{j}(2p+1)!!(2j-2p-1)!! }
{(2j+2)!!(2p+1)}=\frac{1}{4k}-\frac{1}{4}$$

{\bf Theorem 3}.\qquad For any  $k\in N,\,k\geq 1$,
$$
\begin{array}{ll}
\displaystyle{
\sum_{j=0}^{k-1}\sum_{r=0}
^{j}\sum_{p=0}^{j-r}\sum_{s=0}^{r}\sum_{h=0}^{p+s}}&\left(\matrix 2k\\
k+1+j\\ \endmatrix\right)
\small{\left(\matrix j-r\\ p\\ \endmatrix\right)}\left(\matrix r\\ s\\
\endmatrix\right)\left(\matrix p+s\\ h\\\endmatrix\right) \\
&\times\dfrac{(-1)^{j+p+r+1} \Gamma(k-j) \Gamma(1-h+(p+s+j+1)/2) }
{2^{2k-2+p+s-j}\Gamma(k-h+(p+s-j+3)/2)}
=\dfrac{1}{2k+1}-1\end{array}
$$\\

{\bf Acknowledgement}: The author is very grateful to the referee for 
valuable suggestions and corrections which apparently improved this paper.
\\\\
\centerline{\Large References}
\baselineskip 12pt
\begin{description}
\item{[BD]\,} Bremermann. H. J. $\&$ Durand. L., On analytic continuation,
multiplication and Fourier transformations of Schwarts distributions,
J. Math. Phys., 1961, 2, 240-258.
\item{[DR]} F.Diener \& G.Reeb, Analyse Non Standard, Hermann, Editers Sci.
et Arts, 1989
\item{[It1]} Itano, M., On the theory of the multiplicative product of
distributions, J. Sci. Hiroshima Uni., Ser., A-1, 1966, 30, 151-181.
\item{[It2]} Itano, M., On the multiplicative products of $x_{+}^{\alpha}$
and $x_{+}^{\beta}$, J. Sci. Hiroshima Uni., Ser., A-1, 1965, 225-241.
\item{[LK]} Li Ya-Qing $\&$ Kuribayashi,Y., Products $\delta (x_{1})\otimes
I(x_{2},\cdots,x_{n})$ and $\delta (x_{1},\cdots,x_{n})$, J. Fac. Educ. Tottori
Univ. (Nat. Sci.), 42, (1994) 109-117.
\item{[L]} Li Bang-He, Nonstandard analysis and multiplication of
distributions, Scientia Sinica, 1978, 5, 561-585.
\item{[LL1]} Li Bang-He $\&$ Li Ya-Qing, On multiplications of distributions,
Acta scientiarum naturalium universitatis Jilinenesis, 1981, 13-30.
\item{[LL2]} Li Bang-He $\&$ Li Ya-Qing, Nonstandard analysis and multiplication
of distributions in any dimension, Scientia Sinica, 1985, 28, 716-726.
\item{[LL3]} Li Bang-He $\&$ Li Ya-Qing, Multiplication of distributions via 
harmonic representations, Functional Analysis in China, Kluwer Academic 
Publishers, 1996, 90-106.
\item{[LL4]} Li Bang-He $\&$ Li Ya-Qing, On the harmonic and analytic
representations of distributions, Scientia Sinica, 1985, 28, 923-937.
\item{[O]} M.Oberguggenberger, Multiplication of distributions and
applications to partial differential equations, Longman Sci. Technical, UK, 259,
1992
\item{[R]} Robinson, A., Nonstandad analysis, North Hollaid, Amsterdam. 1966.
\end{description}
\noindent Author's address

Institute of Systems Science,
Academia Sinica,
Beijing 100080,
P. R. China  \\ e-mail:
yli @ iss06.iss.ac.cn
\end{document}